\documentclass[11pt,a4paper]{amsart}
\usepackage[usenames]{color}
\usepackage{amsmath, amssymb, amsthm, setspace, mathtools, graphicx, enumitem, tikz-cd, microtype, mathrsfs}
\usepackage{thmtools} % produces remreset obsolete warning
\usepackage[pagebackref,colorlinks,linkcolor=blue,citecolor=blue,urlcolor=blue,hypertexnames=true]{hyperref}
\usepackage{makecell}
\usepackage{longtable}

\declaretheorem[style=plain,name=Theorem]{theorem}
\declaretheorem[style=plain,numberlike=theorem,name=Lemma]{lemma}

\numberwithin{equation}{section}

\newcommand{\M}{\mathrm{M}}                         % matrices
\newcommand{\N}{\mathbb{N}}                         % natural numbers
\newcommand{\Z}{\mathbb{Z}}                         % integers
\newcommand{\Q}{\mathbb{Q}}                         % rationals
\newcommand{\Sp}{\mathrm{Sp}}                       % symplectic group
\newcommand{\one}{\mathbf{1}}                       % neutral element of a group

\DeclareMathOperator{\Span}{Span}

\allowdisplaybreaks

\begin{document}
\date{\today}
\title[Arithmetic monodromy in $\mathrm{Sp}(2n)$]{Arithmetic monodromy in $\mathrm{Sp}(2n)$}
\author{Jitendra Bajpai, Daniele Dona, Martin Nitsche}
\address{Max-Planck-Institut f\"ur Mathematik, Bonn, Germany}
\email{jitendra@mpim-bonn.mpg.de}
\address{Einstein Institute of Mathematics, Hebrew University of Jerusalem, Israel}
\email{daniele.dona@mail.huji.ac.il}
\address{Institute of Algebra and Geometry, Karlsruhe Institute of Technology, Germany}
\email{martin.nitsche@kit.edu}
\subjclass[2010]{Primary: 22E40;  Secondary: 32S40;  33C80}  
\keywords{Hypergeometric group, monodromy representation, symplectic group}

%---------------------------------------------------------------------------
%---------------- Abstract-----------------------------------------------
%---------------------------------------------------------------------------

\begin{abstract}
Based on a result of Singh--Venkataramana, Bajpai--Dona--Singh--Singh gave a criterion for a discrete Zariski-dense subgroup of $\Sp(2n,\Z)$ to be a lattice.
We adapt this criterion so that it can be used in some situations that were previously excluded. We apply the adapted method to subgroups of $\Sp(6,\Z)$ and $\Sp(4,\Z)$ that arise as the monodromy groups of hypergeometric differential equations.
In particular, we show that out of the $40$ maximally unipotent $\Sp(6)$ hypergeometric groups more than half are arithmetic, answering a question of Katz in the negative.
\end{abstract}

\maketitle

%\tableofcontents
%----------------------------------------------------------------------------------------
%	Section - Introduction
%----------------------------------------------------------------------------------------

\section{Introduction}\label{se:intro}

In this article we study groups $\Gamma<\mathrm{Sp}_{\Omega}(\mathbb{Z})$ Zariski-dense in $\mathrm{Sp}_{\Omega}$, where $\Omega$ is a non-degenerate symplectic form on $\Q^{2n}$, $n\geq 2$,  integral on $\Z^{2n}$.
Specifically, our main focus is on the groups $\Gamma=\Gamma(f,g)$ generated by the
companion matrices of polynomials $f,g$ of the form $f=\prod_{j=1}^{2n}(x-e^{2\pi i\alpha_{j}})$ and $g=\prod_{k=1}^{2n}(x-e^{2\pi i\beta_{k}})$, where $\alpha_{j},\beta_{k}\in\mathbb{C}$ with $\alpha_{j}-\beta_{k}\notin\mathbb{Z}$ for all $1\leq j,k\leq 2n$.

These groups arise as the monodromy groups of the hypergeometric differential equations $_{2n}F_{2n-1}$ on the thrice-punctured Riemann sphere, and have been studied for a long time.
For details about the construction and its significance, see \cite{Levelt, Sa14, SV, BDN21}.
We call $\Gamma(f,g)$ the \textit{hypergeometric group} associated to the $2n$-tuples $\alpha=(\alpha_{1},\ldots,\alpha_{2n}),\beta=(\beta_{1},\ldots,\beta_{2n})$.

By the classification in~\cite[Thm.~6.5]{BH}, we know how to determine from $\alpha,\beta$ when the Zariski closure of $\Gamma$ is a symplectic group. Up to conjugation, for any given $n$ there are only finitely many possibilities for $\Gamma$ integral and Zariski-dense in $\mathrm{Sp}(2n)$, and we can produce complete lists of $\alpha,\beta$ yielding such $\Gamma$. Particularly significant are the cases of \textit{maximally unipotent monodromy}, i.e.\ when $\alpha=(0,\ldots,0)$, since these are closely related to families of Calabi--Yau manifolds studied in physics.
For simplicity we also work up to primitivity (see~\cite[Def.~5.1]{BH}) and scalar shift (see~\cite[Def.~5.5]{BH}): under these assumptions, we have
for $n\in\{2,3\}$
\begin{itemize}
\item $58$ symplectic hypergeometric groups in $\mathrm{Sp}(4)$, of which $14$ are maximally unipotent, listed in~\cite[Sect.~5]{SV} (not up to scalar shift) and in~\cite[Sect.~2]{BDN21}, and
\item $458$ symplectic hypergeometric groups in $\mathrm{Sp}(6)$, of which $40$ are maximally unipotent, listed in~\cite[Sect.~4-5-6]{bdss} and~\cite[Sect.~6]{BDSSa}.
\end{itemize}

An important question is to determine whether a given $\Gamma$ is \textit{arithmetic} or \textit{thin}, namely whether it has respectively finite or infinite index inside $\mathrm{Sp}_{\Omega}(\mathbb{Z})$. Answering one way or another has disparate consequences, many of them following from the study of superstrong approximation (see \cite[Sect.~2]{Sa14}).
For $\mathrm{Sp}(4)$, we know from~\cite{SV,S15S} that $7$ of the $14$ hypergeometric groups associated to $\alpha=(0,0,0,0)$ are arithmetic, and the other $7$ are thin by~\cite{BT}. During the workshop on ``Thin Groups and Super Approximation'' held at IAS Princeton in March 2016, N.~Katz raised the question of whether the same symmetry occurs for $\mathrm{Sp}(6)$ as well.

The techniques involved in proving arithmeticity or thinness are to date quite different from each other. For results on thinness, see~\cite{BT, BDN21}. A first sufficient criterion for $\Gamma$ to be arithmetic was proved by Singh and Venkataramana~\cite[Thm.~1.2]{SV}. Later, a second one was derived from the first in~\cite[Prop.~1]{bdss}, which was in practice easier to check. Our first result is another criterion derived from the one of Singh and Venkataramana.

\begin{lemma}\label{lem:main-lemma}
For $n\geq 2$ let $\Omega$ be a non-degenerate symplectic form on $\Q^{2n}$ that is integral on $\Z^{2n}\subset\Q^{2n}$, and let $\Gamma<\Sp_\Omega(\Z)$ be a Zariski-dense subgroup. The following are equivalent:
\begin{enumerate}
\renewcommand{\labelenumi}{\textup{(\theenumi)}}
\item
$\Gamma$ has finite index in $\Sp_\Omega(\Z)$, hence is a lattice in $\Sp_\Omega(\Q)$;
\item
$\Gamma$ contains two transvections $X_i=\one+\lambda_i x_i\Omega(x_i,\cdot)$ such that $x_1,x_2$ are linearly independent and $\Omega$-orthogonal \textup{(}i.e., $\langle X_1,X_2\rangle\cong\Z^2$\textup{)}.
\end{enumerate}
\end{lemma}

We apply Lemma~\ref{lem:main-lemma} in the cases $n=2,3$. Let us start with $\mathrm{Sp}(6)$. Of the $40$ groups associated to $\alpha=(0,0,0,0,0,0)$, we know from \cite{bdss} that at least $18$ are arithmetic. The present article proves that another $5$ groups are arithmetic, settling Katz's question in the negative.

\begin{theorem}\label{thm:sp6-mum}
The hypergeometric groups in $\mathrm{Sp}(6)$ with maximally unipotent monodromy listed in Table~\ref{tab:sp6-mum} are arithmetic.
\end{theorem}

{\centering
\renewcommand{\arraystretch}{1.7}
\begin{longtable}{|c|c|c|}
\caption{Arithmetic groups in $\mathrm{Sp}(6)$ with $\alpha=(0,0,0,0,0,0)$, \\ labelled according to the numeration in \cite[Table~A]{bdss}.}\label{tab:sp6-mum}\\
\hline
Label & $\beta$ & successful word $\gamma$ \\
\hline
\hline
A-15 & $\left(\frac{1}{3},\frac{1}{3},\frac{1}{3},\frac{2}{3},\frac{2}{3},\frac{2}{3}\right)$ & $A^{2}B^{-5}$ \\
\hline
A-16 & $\left(\frac{1}{3},\frac{1}{3},\frac{2}{3},\frac{2}{3},\frac{1}{4},\frac{3}{4}\right)$ & \makecell[c]{\\[-0.4cm]
$ABA^{-1}B^{-8}A\left(AB^{-1}\right)^7A^{-1}B^3\left(AB^{-1}\right)^3\ldots$ \\ $\ldots ABA\left(AB^{-1}\right)^3B^{-1}A\left(BA^{-1}\right)^3BA^2B^4A^5$} \\
\hline
A-21 & $\left(\frac{1}{3},\frac{2}{3},\frac{1}{5},\frac{2}{5},\frac{3}{5},\frac{4}{5}\right)$ & $B^2A^{-1}B^5A^{-1}B\left(BA^{-1}\right)^9A^{-1}B^{-1}AB^{-6}A^4B^{-6}A^2$ \\
\hline
A-24 & $\left(\frac{1}{3},\frac{2}{3},\frac{1}{12},\frac{5}{12},\frac{7}{12},\frac{11}{12}\right)$ & $B^6$ \\
\hline
A-39 & $\left(\frac{1}{14},\frac{3}{14},\frac{5}{14},\frac{9}{14},\frac{11}{14},\frac{13}{14}\right)$ & $ABA^{-1}B^{-3}A^{-4}B^{-3}$ \\
\hline
\end{longtable}}

We also show in another article~\cite{BDN21} that the remaining $17$ groups are thin, completing the classification for maximally unipotent groups in $\mathrm{Sp}(6)$.

Then, we address the question of arithmeticity among the remaining $418$ symplectic hypergeometric groups with $\alpha\neq(0,0,0,0,0,0)$. We know from~\cite{SV,DFH,bdss,JB} that at least $361$ of them are arithmetic. We can now show the arithmeticity of $8$ more. Since the proof in~\cite{JB} (for the groups C-01, C-10, C-42, C-59 in Table~\ref{tab:sp6-other}) is rather long and technical, we also reprove arithmeticity for these groups with our new criterion.

\begin{theorem}\label{thm:sp6-other}
The hypergeometric groups in $\mathrm{Sp}(6)$ listed in Table~\ref{tab:sp6-other} are arithmetic.
\end{theorem}

{\centering
\renewcommand{\arraystretch}{1.7}
\begin{longtable}{|c|c|c|c|}
\caption{Arithmetic groups in $\mathrm{Sp}(6)$ with $\alpha\neq(0,0,0,0,0,0)$, \\ labelled according to the numeration in \cite[Table~C]{bdss}.}\label{tab:sp6-other}\\
\hline
Label & $\alpha$ & $\beta$ & successful word $\gamma$ \\
\hline
\hline
C-01 & $\left(0,0,0,0,\frac{1}{2},\frac{1}{2}\right)$ & $\left(\frac{1}{3},\frac{1}{3},\frac{2}{3},\frac{2}{3},\frac{1}{6},\frac{5}{6}\right)$ & $BA$ \\
\hline
C-09 & $\left(0,0,0,0,\frac{1}{3},\frac{2}{3}\right)$ & $\left(\frac{1}{7},\frac{2}{7},\frac{3}{7},\frac{4}{7},\frac{5}{7},\frac{6}{7}\right)$ & $BA^{-1}B^{-1}AB^{-2}AB^{-1}\left(B^3A^2\right)^4$ \\
\hline
C-10 & $\left(0,0,0,0,\frac{1}{3},\frac{2}{3}\right)$ & $\left(\frac{1}{9},\frac{2}{9},\frac{4}{9},\frac{5}{9},\frac{7}{9},\frac{8}{9}\right)$ & $BAB^2 A^{-4}$ \\
\hline
C-29 & $\left(0,0,0,0,\frac{1}{6},\frac{5}{6}\right)$ & $\left(\frac{1}{3},\frac{1}{3},\frac{1}{3},\frac{2}{3},\frac{2}{3},\frac{2}{3}\right)$ & $B^4A^{-1}B^6A^{-1}BA^{-2}$ \\
\hline
C-30 & $\left(0,0,0,0,\frac{1}{6},\frac{5}{6}\right)$ & $\left(\frac{1}{3},\frac{1}{3},\frac{2}{3},\frac{2}{3},\frac{1}{4},\frac{3}{4}\right)$ & $B^3A^7BA^{-1}BA^{-1}B^4AB^4A^{-1}BA^2$ \\
\hline
C-31 & $\left(0,0,0,0,\frac{1}{6},\frac{5}{6}\right)$ & $\left(\frac{1}{3},\frac{2}{3},\frac{1}{5},\frac{2}{5},\frac{3}{5},\frac{4}{5}\right)$ & $BABA^{-3}\left(BA^{-1}\right)^3B^{-3}A^{-2}BA^{-1}B^6A^{-7}$ \\
\hline
C-39 & $\left(0,0,\frac{1}{3},\frac{2}{3},\frac{1}{6},\frac{5}{6}\right)$ & $\left(\frac{1}{7},\frac{2}{7},\frac{3}{7},\frac{4}{7},\frac{5}{7},\frac{6}{7}\right)$ & $B^2A^{-2}B^{-3}A^{-2}B^{-3}A^2B^2A$ \\
\hline
C-42 & $\left(0,0,\frac{1}{4},\frac{1}{4},\frac{3}{4},\frac{3}{4}\right)$ & $\left(\frac{1}{3},\frac{2}{3},\frac{1}{12},\frac{5}{12},\frac{7}{12},\frac{1}{12}\right)$ & $A^3$ \\
\hline
C-51 & $\left(0,0,\frac{1}{6},\frac{1}{6},\frac{5}{6},\frac{5}{6}\right)$ & $\left(\frac{1}{2},\frac{1}{2},\frac{1}{12},\frac{5}{12},\frac{7}{12},\frac{11}{12}\right)$ & $A^6$ \\
\hline
C-59 & $\left(0,0,\frac{1}{12},\frac{5}{12},\frac{7}{12},\frac{11}{12}\right)$ & $\left(\frac{1}{3},\frac{2}{3},\frac{1}{4},\frac{1}{4},\frac{3}{4},\frac{3}{4}\right)$ & $BA$ \\
\hline
C-60 & $\left(\frac{1}{3},\frac{1}{3},\frac{1}{3},\frac{2}{3},\frac{2}{3},\frac{2}{3}\right)$ & $\left(\frac{1}{6},\frac{1}{6},\frac{1}{6},\frac{5}{6},\frac{5}{6},\frac{5}{6}\right)$ & $BA^{-1}BA^2BAB^{-1}A^{-4}$ \\
\hline
C-61 & $\left(\frac{1}{3},\frac{1}{3},\frac{1}{3},\frac{2}{3},\frac{2}{3},\frac{2}{3}\right)$ & $\left(\frac{1}{9},\frac{2}{9},\frac{4}{9},\frac{5}{9},\frac{7}{9},\frac{8}{9}\right)$ & $AB^{-4}A^3$ \\
\hline
\end{longtable}}

Together with the thinness results in~\cite{BDN21}, there are only $3$ groups left in $\mathrm{Sp}(6)$ whose status as arithmetic or thin is still unknown.

Finally, we turn to $\mathrm{Sp}(4)$. Of the $58$ symplectic hypergeometric groups in $\mathrm{Sp}(4)$, $53$ are already known to be either arithmetic or thin by~\cite{SV, BT, S15S, S17, BSS}: our article shows the arithmeticity of $2$ more.

\begin{theorem}\label{thm:sp4}
The hypergeometric groups in $\mathrm{Sp}(4)$ listed in Table~\ref{tab:sp4} are arithmetic.
\end{theorem}

{\centering
\renewcommand{\arraystretch}{1.7}
\begin{longtable}{|c|c|c|c|}
\caption{Arithmetic groups in $\mathrm{Sp}(4)$ labelled according to the numeration in~\cite[Table~4]{BDN21}.}\label{tab:sp4}\\%, labelled following~\cite{SV}*{Table~2}
\hline
Label & $\alpha$ & $\beta$ & successful word $\gamma$ \\
\hline
\hline
30 & $\left(0,0,\frac{1}{4},\frac{3}{4}\right)$ & $\left(\frac{1}{5},\frac{2}{5},\frac{3}{5},\frac{4}{5}\right)$ & $BA^2B^{-2}(A^{-2}B^{-2}A^3B^{-2})^2$ \\
\hline
40 & $\left(0,0,\frac{1}{6},\frac{5}{6}\right)$ & $\left(\frac{1}{8},\frac{3}{8},\frac{5}{8},\frac{7}{8}\right)$ & $(AB^{-1})^2A^{-1}(BA)^{-2}A^{-3}B^{-2}A^4B^{-2}A^{-4}B^{-2}A^4B^{-2}A^{-3}$ \\
\hline
\end{longtable}}

In~\cite{BDN21} the remaining $3$ groups are proved to be thin, completing the classification of the $\mathrm{Sp}(4)$ hypergeometric groups.

\section{Criterion for proving arithmeticity}
Recall that a transvection is a matrix $X$ such that $X-\one$ has rank one. When the symplectic form $\Omega$ is non-degenerate, any transvection $X\in\Sp_\Omega$ has the form $X=\one+\lambda x\Omega(x,\cdot)$ for some $\lambda\in\mathbb{R}$ and $x\in\mathbb{R}^{2n}$. The arithmeticity criterion of Singh--Venkataramana mentioned in the introduction is as follows.

\begin{theorem}[Singh--Venkataramana~\cite{SV}]\label{thm:sv}
Assume that $\Gamma$ contains three transvections $X_i=\one+\lambda_i x_i\Omega(x_i,\cdot)$, 
whose directions span a three-dimensional subspace $W=\Span(\{x_i\})\subset\Q^n$, such that the restriction $\Omega_{|W}$ of the symplectic form is non-trivial. Assume further that the group generated by the restrictions ${X_i}_{|W}$ contains a non-trivial element of the unipotent radical of $\Sp_{\Omega_{|W}}(\Q)$.

Then $\Gamma$ has finite index in $\Sp_\Omega(\Z)$, hence is a lattice in $\Sp_\Omega(\Q)$.
\end{theorem}

In~\cite{bdss}, Bajpai, Dona, Singh and Singh derived from this result a more specialized criterion that is easier to check and implies that the conditions of Theorem~\ref{thm:sv} are met. They used it to develop a computer algorithm that takes a transvection $X_1\in\Gamma$ and searches for two $\Gamma$-conjugates $X_2,X_3\in\Gamma$ that satisfy the conditions of Theorem~\ref{thm:sv}. This algorithm is successful in many cases, but it can only be used when the entries of $X_1-\one$ have greatest common divisor $\leq 2$. The advantage of our new criterion is that it does not suffer from this restriction.

\newtheorem*{restate:main-lemma}{Lemma~\ref{lem:main-lemma}}
\begin{restate:main-lemma}
For $n\geq 2$ let $\Omega$ be a non-degenerate symplectic form on $\Q^{2n}$ that is integral on $\Z^{2n}\subset\Q^{2n}$, and let $\Gamma<\Sp_\Omega(\Z)$ be a Zariski-dense subgroup. The following are equivalent:
\begin{enumerate}
\renewcommand{\labelenumi}{\textup{(\theenumi)}}
\item
$\Gamma$ has finite index in $\Sp_\Omega(\Z)$, hence is a lattice in $\Sp_\Omega(\Q)$;
\item
$\Gamma$ contains two transvections $X_i=\one+\lambda_i x_i\Omega(x_i,\cdot)$, such that $x_1,x_2$ are linearly independent and $\Omega$-orthogonal \textup{(}i.e., $\langle X_1,X_2\rangle\cong\Z^2$\textup{)}.
\end{enumerate}
\end{restate:main-lemma}

\begin{proof}
\textbf{1) $\Rightarrow$ 2)}
We choose two linearly independent and $\Omega$-orthogonal elements $x_i\in\mathbb{Z}^{2n}$. The two transvections $X_i=\one+x_i\Omega(x_i,\cdot)$ lie in $\Sp_\Omega(\Z)$, and since $\Gamma<\Sp_\Omega(\Z)$ has finite index, each of them has a power that lies in $\Gamma$.

\textbf{2) $\Rightarrow$ 1)}
First, consider the set $\{\gamma\in\Sp(\Q)\mid\gamma x_1\perp x_1\vee\gamma x_1\perp x_2\}$. This is a non-trivial Zariski-closed subset of $\Sp(\Q)$. We pick any element $\gamma$ in the intersection of its complement with the Zariski-dense subset $\Gamma$. Then $X_3\vcentcolon=\gamma X_1\gamma^{-1}\in\Gamma$ is a transvection whose direction $x_3$ is neither orthogonal to $x_1$ nor to $x_2$.
The directions $x_1,x_2,x_3$ span a $3$-dimensional subspace $W$ such that $\Omega_{|W}$ is non-trivial.

Next, we define $R\in\langle X_1,X_2\rangle <\Gamma$ by $R\vcentcolon={X_1}^{N\lambda_2\Omega(x_2,x_3)^2}{X_2}^{-N\lambda_1\Omega(x_1,x_3)^2}$, where $N\in\N$ is chosen large enough to make the exponents integers.
$R$ acts as the identity on $x_1,x_2$, and it sends $x_3$ to $x_3+N\lambda_1\lambda_2\Omega(x_1,x_3)\Omega(x_2,x_3)^2x_1-N\lambda_1\lambda_2\Omega(x_1,x_3)^2\Omega(x_2,x_3)x_2$. Hence, the restriction $R_{|W}$ is a transvection in the direction $\Omega(x_2,x_3)x_1-\Omega(x_1,x_3)x_2$, which is contained in the subspace $W\cap W^\perp$. Moreover, no other direction can be in $W\cap W^\perp$ because of our choice of $x_{3}$, so $W\cap W^\perp$ is $1$-dimensional.

The subgroup $S<\Sp_{\Omega_{|W}}(\Q)$ consisting of all transvections with direction $W\cap W^\perp$ is normal in $\Sp_{\Omega_{|W}}(\Q)$ by construction, and it is also abelian because $\dim(W\cap W^\perp)=1$. Hence, $S$ lies in the radical of $\Sp_{\Omega_{|W}}(\Q)$. Therefore, $R_{|W}$ is a non-trivial element of the unipotent radical of $\Sp_{\Omega_{|W}}(\Q)$ and all conditions of Theorem~\ref{thm:sv} are satisfied.
\end{proof}

\section{Application of the arithmeticity criterion}
Finally, we explain how to apply Lemma~\ref{lem:main-lemma} to the groups whose parameters are listed in Tables~\ref{tab:sp6-mum}, \ref{tab:sp6-other}, \ref{tab:sp4}.

The hypergeometric groups $\Gamma<\Sp(2n)$ come with two generating matrices $A,B$ such that $T\vcentcolon=A^{-1}B$ is a transvection. We can write $T=\one+v_R v_L$ for $v_R\in\M_{2n\times 1}(\Z)$ and $v_L\in\M_{1\times 2n}(\Z)$. With the computer we search for group elements $\gamma$, written as words in $A^{\pm 1},B^{\pm 1}$, such that $T,\gamma T\gamma^{-1}$ satisfy the condition of Lemma~\ref{lem:main-lemma}. This means that ${v_L}\gamma v_R=0$ and that $v_R,\gamma v_R$ are linearly independent.

The search is done by computing iteratively the sequence of sets $R_0\vcentcolon=\{v_R\}$, $R_{i+1}\vcentcolon=R_i\cup AR_i\cup A^{-1}R_i\cup BR_i\cup B^{-1}R_i$ and checking the two conditions for all their members. To speed up the search, we actually compute $R_{i+1}=R_i\cup A(R_i\setminus R_{i-1})\cup A^{-1}(R_i\setminus R_{i-1})\cup B(R_i\setminus R_{i-1})\cup B^{-1}(R_i\setminus R_{i-1})$. We also filter out all vectors whose entries exceed a fixed bound in absolute value. This last heuristic drastically prunes the search tree, allowing for deeper searches and guaranteeing termination of the search. Furthermore, it averts the risk of rounding errors resulting from numerical computations with too large numbers.

A MATLAB-implementation of the search algorithm can be found on the third-named author's GitHub page\footnote{\url{https://github.com/MartinNitsche/Arithmeticity-in-Sp}}. Each computation took less than a minute on standard consumer hardware. The words $\gamma$ found through this method are listed next to each entry in Tables~\ref{tab:sp6-mum}, \ref{tab:sp6-other}, \ref{tab:sp4}.

It is elementary to check that for these $\gamma$ the transvections $T,\gamma T\gamma^{-1}$ satisfy the preconditions of Lemma~\ref{lem:main-lemma}, and hence Theorems~\ref{thm:sp6-mum},~\ref{thm:sp6-other},~\ref{thm:sp4} follow.

%----------------------------------------------------------------------------------------
%	Acknowledgements
%----------------------------------------------------------------------------------------
\section*{Acknowledgements}

JB is supported through the fellowship from MPIM, Bonn. DD is supported by the Israel Science Foundation Grants No. 686/17 and 700/21 of A.~Shalev, and the Emily Erskine Endowment Fund; he has been a postdoc at the Hebrew University of Jerusalem under A.~Shalev in 2020/21 and 2021/22. MN was supported by ERC Consolidator grant 681207 and DFG grant 281869850 (RTG 2229, ``Asymptotic Invariants and Limits of Groups and Spaces'').

%----------------------------------------------------------------------------------------
%                  Bibliogrpahy
%----------------------------------------------------------------------------------------
\nocite{}
\bibliographystyle{abbrv}
\bibliography{BDN}

\end{document}